\documentclass [11pt]{article}
\usepackage {amsfonts}
\usepackage {mathrsfs}
\usepackage {amsmath}
\usepackage {latexsym}
\usepackage {amssymb}
\usepackage {amsthm}
\usepackage {amscd}
\date{}
\renewcommand{\baselinestretch}{1.2}
\title{\bf Classification of quasi-affine Generalized Dynkin Diagrams  with Rank $3$ and Rank $2$}
\author{ \small Zhengtang Tan $^{a}$,  Shouchuan Zhang $^{b}$   \\
\small $a$.School of Engineering and Design, Hunan Normal University\\ Changsha  410081,   P.R. China \\
\small $b$. Department  of Mathematics,   Hunan University\\
Changsha  410082,  P.R. China\\
\small {\tt Emails:  z9491@sina.cn (SZ);   1843186255@qq.com (ZTT)} }
\date{}

\begin{document}
\newtheorem{Proposition}{Proposition}[section]
\newtheorem{Theorem}[Proposition]{Theorem}
\newtheorem{Definition}[Proposition]{Definition}
\newtheorem{Corollary}[Proposition]{Corollary}
\newtheorem{Lemma}[Proposition]{Lemma}
\newtheorem{Example}[Proposition]{Example}
\newtheorem{Remark}[Proposition]{Remark}

\maketitle


\begin {abstract}  All quasi-affine  connected Generalized Dynkin Diagram with rank $= 3$ and $2$ are found. All quasi-affine  Nichols (Lie braided) algebras with rank $ 3$ and $2$ are also found.

\end {abstract}

\section {Introduction and Preliminaries}\label {s0}

Nichols algebras play a fundamental role in the classification of finite-dimensional complex pointed Hopf algebras by means of the lifting method developed by Andruskiewitsch and Schneider \cite {AS02, AS10, AHS08}.
Heckenberger \cite {He06a, He05} classified arithmetic root systems.   Heckenberger \cite {He06b} proved a GDD is  arithmetic  if and only if corresponding matrix is a finite Cartan matrix for GDDs of  Cartan types. W. Wu,   S. Zhang and   Y.-Z. Zhang \cite {WZZ15b} proved a Nichols Lie braided algebra is a finite dimensional  if and only if its GDD, which   fixed
parameter is of finite order, is  arithmetic.


We now recall some basic concepts of the graph theory (see \cite {Ha}).
Let $\Gamma _1$ be a non-empty set and $\Gamma _2 \subseteq \{  \{ u, v\} \mid u, v \in \Gamma _1, \hbox { with } u \not= v \} \subseteq 2 ^{\Gamma_1}.$ Then $\Gamma = (\Gamma _1, \Gamma_2)$ is called a graph;  $\Gamma_1$ is called the vertex set of  $\Gamma$; $\Gamma_2$ is called the edge set of  $\Gamma$; Element $\{u, v\} \in \Gamma_2$ is called an edge. Let $F$ be an algebraically closed field of characteristic zero and $F^* := \{x \mid x\in F, x \not= 0 \}$.  If   $\{x_1,   \cdots,   x_n\}$ is  a basis of  vector space  $V$ and
$C(x_i\otimes  x_j) = q_{ij} x_j\otimes x_i$ with $q_{ij} \in F^*$,
then $V$  is called a braided vector space of diagonal type, $\{x_1,   \cdots,   x_n\}$
is called  canonical basis and $(q_{ij})_{n\times n}$  is called braided matrix. Let
$\widetilde{q}_{i j}:= q_{ij}q_{ji}$ for  $i, j \in \{1, 2, \cdots, n \}$ with $i \not=j$. If let $\Gamma _1= \{1, 2, \cdots, n \}$ and $\Gamma _2 = \{  \{ u, v\} \mid \widetilde{q}_{u v} \not= 1,   u, v \in \Gamma _1, \hbox { with } u \not= v \}$, then $\Gamma = (\Gamma _1, \Gamma_2)$ is a graph. Set $q_{ii}$ over vertex $i$ and $\widetilde{q}_{i j}$ over edge $a_{ij}$ for $i, j \in \Gamma _1$ with $i \not=j$ when $\widetilde{q}_{i j} \not=1$.
Then $\Gamma = (\Gamma _1, \Gamma_2)$ is called a Generalized Dynkin Diagram of braided vector space $V,$ written a GDD in short( see \cite [Def. 1.2.1] {He05}).
If $\Delta (\mathfrak{B}(V))$ is an arithmetic root system,   then we call its GDD  an arithmetic GDD  for convenience.

Let $(q_{ij})_{n \times n }$ be a braided  matrix. If  $q_{ij} q_{ji}  \left \{  \begin{array}{ll}
 \not= 1,  & \mbox {when }  \mid j-i \mid = 1\\
  =1,   & \mbox {when}  \mid j-i \mid \not= 1    \\
\end{array}\right. $ for any $1\le i\not= j \le n,$ then  $(q_{ij})_{n \times n }$ is called a  chain or labelled  chain.
If $(q_{ij})_{n \times n }$ is  a  chain and \begin {eqnarray} \label {ppe1}(q_{11} q_{1, 2} q_{2, 1} -1)(q_{11}  + 1)=0;
 (q_{n, n } q_{n, n-1} q_{n-1, n} -1)(q_{n,n}  + 1)=0;\end {eqnarray} i.e.
  \begin {eqnarray} \label {ppe2}
   q_{ii}  + 1= q_{i, i-1} q_{i-1, i}q_{i, i + 1} q_{i + 1, i}-1=0
   \end  {eqnarray}
    \begin {eqnarray} \label {ppe3}
   \mbox { or  }
  q_{ii}q_{i, i-1} q_{i-1, i}=q_{ii} q_{i, i + 1} q_{i + 1, i}=1,\end  {eqnarray}  $1<i < n$,
   then  the braided matrix's GDD is called a simple chain (see \cite [Def.1]{He06a}).  Conditions (\ref {ppe1}), (\ref {ppe2}) and (\ref {ppe3}) are called simple chain conditions.  Let
$q:= q_{n, n} ^2 q_{n -1, n}q_{n, n-1}$.

 For $0\le j \le n$ and $0<i_1 < i_2 < \cdots < i_j \le n,$ let $C_{n,q, i_1, i_2, \cdots, i_j }$ denote  a simple chain which satisfies the condition: $\widetilde{q}_{i, i-1}=q$ if and only if  $i \in \{i_1, i_2, \cdots, i_j\}$,
where
$\widetilde{q} _{1, 0}:= \frac {1}{q_{11}^2\widetilde{q}_{12}}$. For example,
 $n \ge 2$,

$\begin{picture}(100,      15)
\put(27,      1){\makebox(0,     0)[t]{$\bullet$}}
\put(60,      1){\makebox(0,      0)[t]{$\bullet$}}
\put(93,     1){\makebox(0,     0)[t]{$\bullet$}}
\put(159,      1){\makebox(0,      0)[t]{$\bullet$}}
\put(192,     1){\makebox(0,      0)[t]{$\bullet$}}
\put(28,      -1){\line(1,      0){30}}
\put(61,      -1){\line(1,      0){30}}
\put(130,     1){\makebox(0,     0)[t]{$\cdots\cdots\cdots\cdots$}}
\put(160,     -1){\line(1,      0){30}}
\put(22,     -15){1}
\put(58,      -15){2}
\put(91,      -15){3}
\put(157,      -15){$n-1$}
\put(191,      -15){$n$}
\put(22,     10){$q$}
\put(58,      10){$q$}
\put(91,      10){$q$}
\put(157,      10){$q$}
\put(191,      10){$q$}
\put(40,      5){$q^{-1}$}
\put(73,      5){$q^{-1}$}
\put(172,     5){$q^{-1}$}
\put(210,        -1)  {$, q \in F^{*}/\{1, -1\}$. }
\end{picture}$\\

is  $C_{n,q, i_1, i_2, \cdots, i_j }$ with $j=0$.\\

$\begin{picture}(100,      15)
\put(27,      1){\makebox(0,     0)[t]{$\bullet$}}
\put(60,      1){\makebox(0,      0)[t]{$\bullet$}}
\put(93,     1){\makebox(0,     0)[t]{$\bullet$}}
\put(159,      1){\makebox(0,      0)[t]{$\bullet$}}
\put(192,     1){\makebox(0,      0)[t]{$\bullet$}}
\put(28,      -1){\line(1,      0){30}}
\put(61,      -1){\line(1,      0){30}}
\put(130,     1){\makebox(0,     0)[t]{$\cdots\cdots\cdots\cdots$}}
\put(160,     -1){\line(1,      0){30}}
\put(22,     -15){1}
\put(58,      -15){2}
\put(91,      -15){3}
\put(157,      -15){$n-1$}
\put(191,      -15){$n$}
\put(22,     10){$q^{-1}$}
\put(58,      10){$q^{-1}$}
\put(91,      10){$q^{-1}$}
\put(157,      10){$q^{-1}$}
\put(191,      10){$-1$}
\put(40,      5){$q^{}$}
\put(73,      5){$q^{}$}
\put(172,     5){${q}$}
\put(210,        -1)  {$, q \in F^{*}/\{1\}$. }
\end{picture}$\\

is  $C_{n,q, i_1, i_2, \cdots, i_j }$ with $j=n$.

For the convenience, we let
 $C_{1,q, i_1, i_2, \cdots, i_j }$  denote the GDD with length $1$ satisfied the following conditions:   $q= q_{11}$ when $q_{11} \not=-1$; $q$ can be any number   when $q_{11} =-1$.

Every connected subGDD of
every   arithmetic GDD in Row 1-10 in Table C is  called a classical GDD.  Every connected     GDD which is not classical is called an exception.

If omitting every vertex in a connected GDD with rank $n>1$  is an arithmetic GDD and this GDD is not an arithmetic  GDD, then this GDD is called  a quasi-affine  GDD over  connected subGDD of this GDD with rank $n-1$, quasi-affine  GDD in short. In this case, Nichols algebra $\mathfrak B(V)$
and Nichols Lie braded algebra $\mathfrak L(V)$ are said to be quasi-affine.
In other word, if a GDD is quasi-affine  of a braided vector space $V$ which   fixed
parameter is of finite order, then Nichols algebra and Nichols Lie braded algebra
of every proper subGDD are finite dimensional with $\dim   \mathfrak B { (V)} = \infty$
and $\dim   \mathfrak L { (V)} = \infty$.

In this paper,  using Table $A1$, $A2$, $B$ and $C$ in \cite {He06a, He05}, we find all quasi-affine  connected Generalized Dynkin Diagrams with rank $=3$ and $2$. We also  find all quasi-affine   Nichols  algebras and quasi-affine  Nichols Lie braided algebras with rank $ 3$ and $2$.

In this paper we use following notations without special announcement: every braided vector  $V$  is of diagonal type with dimension $n > 1$; is connected with fixed parameter $q \in R_3$. Let  $F$ be an  algebraic closed base field of characteristic zero. $\mathbb Z =: \{x \mid x \hbox { is an integer }\}.$
$\mathbb N_0 =: \{x \mid x \in \mathbb Z,    x\ge 0\}. $

\section {Properties about arithmetic GDD  } \label {s2}

\begin {Lemma} \label {mainlemma} A GDD is a classical GDD if and only if it is one of classical Type 1-7, Here
classical types are  listed as follows:\\

 {\ }\ \ \ \ \ \ \ \ \ \ \ \ \ \ \ \ \ \ \ \ \ \ \ \   \ \ \ \ \ \  $\begin{picture}(100,      15)

\put(-125,      -1){ {\rm   Type   1},  $2\le n$.}
\put(80,      1){\makebox(0,      0)[t]{$\bullet$}}

\put(48,      -1){\line(1,      0){33}}
\put(10,      1){\makebox(0,     0)[t]{$C_{n-1, q, i_1, i_2, \cdots, i_j }$}}

\put(-18,     10){$$}
\put(0,      5){$$}
\put(22,     10){$$}
\put(50,      5){$q^{-2}$}

\put(68,      10){$q^2$}

 \ \ \ \ \ \ \ \ \ \ \ \ \ \ \ \ \ \ \
  \ \ \ \ \ \ \ \ \ \ \ \ \ \ \ \ \ \ \ { }$q \in F^{*}\setminus \{1, -1\}$, $0\le j\le n-1.$

 \put(220,      -1) {}
\end{picture}$\\

 {\ }\ \ \ \ \ \ \ \ \ \ \ \ \ \ \ \ \ \ \ \ \ \ \ \   \ \ \ \ \ \  $\begin{picture}(100,      15)

\put(-125,      -1){ {\rm   Type   2}. $2\le n$.}

\put(90,      1){\makebox(0,      0)[t]{$\bullet$}}

\put(58,      -1){\line(1,      0){33}}
\put(27,      1){\makebox(0,     0)[t]{$C_{n-1, q^2, i_1, i_2, \cdots, i_j }$}}

\put(-8,     10){$$}
\put(0,      5){$$}
\put(22,     10){$$}
\put(70,      5){$q^{-2}$}

\put(88,      10){$q$}

  \ \ \ \ \ \ \ \ \ \ \ \ \ \ \ \ \ \ \
  \ \ \ \ \ \ \ \ \ \ \ \ \ \ \ \ \ \ \ $q \in F^{*}\setminus \{1, -1\}$.  $0\le j\le n-1.$
\end{picture}$\\

 {\ }\ \ \ \ \ \ \ \ \ \ \ \ \ \ \ \ \ \ \ \ \ \ \ \   \ \ \ \ \ \  $\begin{picture}(100,      15)

\put(-125,      -1){ {\rm   Type   3}, $2\le n$.}

\put(80,      1){\makebox(0,      0)[t]{$\bullet$}}

\put(48,      -1){\line(1,      0){33}}
\put(27,      1){\makebox(0,     0)[t]{$C_{n-1, q^{-2}, i_1, i_2, \cdots, i_j }$}}

\put(-18,     10){$$}
\put(0,      5){$$}
\put(22,     10){$$}
\put(60,      5){$q^{2}$}

\put(78,      10){$-q^{-1}$}

  \ \ \ \ \ \ \ \ \ \ \ \ \ \ \ \ \ \ \
  \ \ \ \ \ \ \ \ \ \ \ \ \ \ \ \ \ \ \ $q \in F^{*}\setminus \{1, -1\}$, $0\le j\le n-1.$

\end{picture}$\\

 {\ }\ \ \ \ \ \ \ \ \ \ \ \ \ \ \ \ \ \ \ \ \ \ \ \   \ \ \ \ \ \  $\begin{picture}(100,      15)

\put(-125,      -1){{\rm   Type   4}. $2\le n$. }

\put(80,      1){\makebox(0,      0)[t]{$\bullet$}}

\put(48,      -1){\line(1,      0){33}}
\put(27,      1){\makebox(0,     0)[t]{$C_{n-1, -q^{-1}, i_1, i_2, \cdots, i_j }$}}

\put(-18,     10){$$}
\put(0,      5){$$}
\put(22,     10){$$}
\put(60,      5){$-q$}

\put(78,      10){$q$}

  \ \ \ \ \ \ \ \ \ \ \ \ \ \ \ \ \ \ \
  \ \ \ \ \ \ \ \ \ \ \ \ \ \ \ \ \ \ \ $q^3 =1$. $0\le j\le n-1.$

\end{picture}$\\
\\

\ \ \ \     \ \ \ \ \ \  $\begin{picture}(100,      15)

\put(-45,      -1){ {\rm   Type   5}. $3\le n$.}
\put(124,      1){\makebox(0,      0)[t]{$C_{n-2, q^{}, i_1, i_2, \cdots, i_j }$}}
\put(190,     -11){\makebox(0,     0)[t]{$\bullet$}}
\put(190,    15){\makebox(0,     0)[t]{$\bullet$}}
\put(162,    -1){\line(2,     1){27}}
\put(190,      -14){\line(-2,     1){27}}

\put(120,      10){$$}

\put(135,      5){$$}

\put(155,     10){$$}

\put(160,     -20){$q^{-1}$}
\put(165,      15){$q^{-1}$}

\put(193,      -12){$q$}
\put(193,      18){$q$}

 \ \ \ \ \ \ \ \ \ \ \ \ \ \ \ \ \ \ \ \ \ \ \ \ \ \ \ \ \ \ \ \ \ \ \ \ \ \
  \ \ \ \ \ \ \ \ \ \ \ \ \ \ \ \ \ \ \ $q\not= 1,$    $0\le j\le n-2.$

\put(215,        -1)  { }
\end{picture}$\\ \\

\ \ \ \      \ \ \ \ \ \  $\begin{picture}(100,      15)

\put(-45,      -1){{\rm   Type   6}. $3\le n$. }

\put(104,      1){\makebox(0,      0)[t]{$C_{n-2, q^{}, i_1, i_2, \cdots, i_j }$}}
\put(170,     -11){\makebox(0,     0)[t]{$\bullet$}}
\put(170,    15){\makebox(0,     0)[t]{$\bullet$}}

\put(142,    -1){\line(2,     1){27}}
\put(170,      -14){\line(-2,     1){27}}

\put(170,      -14){\line(0,     1){27}}

\put(100,      10){$$}
\put(120,      5){$$}

\put(127,     10){$$}

\put(140,     -20){$q^{-1}$}
\put(145,      15){$q^{-1}$}

\put(178,      -20){$-1$}

\put(178,      0){$q^{2}$}

\put(175,      10){$-1$}

 \ \ \ \ \ \ \ \ \ \ \ \ \ \ \ \ \ \ \ \ \ \ \ \ \ \ \ \ \ \ \ \ \ \ \ \ \ \
  \ \ \ \ \ \ \ \ \ \ \ \ \ \ \ \ \ \ \  $q^2 \not= 1$. $0\le j\le n-2.      $
\put(195,        -1)  {.    }
\end{picture}$\\

{\ }\!\!\!\!\!\!\!\!\!\!\!\!\!\!\!\!\!
{\rm   Type   7}. $1\le n$.
$C_{n, q^{-1}, i_1, i_2, \cdots, i_j }.$
$q \not= 1$,  $0\le j\le n$.

\end {Lemma}
{\bf Proof.} We only consider the case $n>4$ since every  GDD classical  type   with $n<5$ is subGDD of a GDD classical type with  $n> 4.$

Necessity is obviously. Now we show the sufficiency.

 Type   1.   Type   1 is classical  by  GDD $1$ of Row $10$ in Table C when  $2 \le j \le n-1$.   Type   1 is classical  by  GDD $1$ of Row $9$ in Table C when  $j=0$.  Type   1 is classical  by GDD $2$ of Row $9$  in Table C and by GDD $1$ of Row $7$  in Table C
 when  $j=1$.

Type 2.  Type 2 is classical  by  GDD $1$ of Row $4$ in Table C when $q^4\not=1$ and   $1 \le j \le n-1$. Type 2 is classical  by  GDD $1$ of Row $3$ in Table C when $q^2\not=1$ and $j=0$. Type 2 is classical  by  GDD $1$ of Row $3$ in Table C when $1 \le j \le n-1.$

Type 3.  Type 3 is classical  by  GDD $2$ of Row $4$ in Table C when $q^4\not=1$ and  $1 \le j \le n-1$.

Type 4.  Type 4 is classical  by  GDD $1$, $2$ of Row $6$ in Table C when  $1 \le j \le n-1$.  Type 4 is classical  by  GDD $1$ of Row $5$ in Table C when  $j=0$.

Type 5.  Type 5 is classical  by  GDD $3$ of Row $10$ in Table C when  $1 \le j \le n-2$ with $q \not=1$.  Type 5 is classical  by  GDD $1$ of Row $8$  in Table C when  $j=0$ with $q^2\not=1$.

Type 6.  Type 6 is classical  by  GDD $2$ of Row $10$ in Table C when  $1 \le j \le n-2$.  Type 1 is classical  by  GDD $3$ of Row $9$ in Table C when  $j=0$. $\Box$

\section {    Quasi-affine  Nichols  (Lie braided) Algebras with Rank $3$    } \label {s2}

In this section all quasi-affine  connected Generalized Dynkin Diagram with rank $= 3 $ are found. All quasi-affine  Nichols (Lie braided) algebras with rank $3 $  are also found.

\begin {Theorem} \label {Mainrank3}
        A connected  GDD with rank $= 3$ is quasi-affine  if and only if it is  one of following lists.


 {\ }

 \  \ \  \ \ \ \ \  \ \ $
$\\ \\ except Case 1:   $\lambda \mu q =1;$   Case 2: $ q^{-1} =\mu ^{2}$, $ \lambda^{-1} =\mu^{-1}$.

\end {Theorem}

{\bf Proof.}
        A connected  GDD with rank 3 is quasi-affine if and only if it is not   arithmetic and   omitting every  vertex still is  arithmetic.

If a GDD is quasi-affine over two GDDs, we write it in below one.

\subsection* {
GDD 1 of Row 6 in Table A1}\label {sub3.2}

{\rm (i)}     Adding on Vertex 2 by a GDD in Table A1.

 We find all quasi-affine GDDs  adding a GDD in Table A1 on Vertex 2 of GDD $1$ of Row $6$. We have to consider all  GDDs adding  a GDD in Table A1 on Vertex 2 of GDD $1$ of Row $6$. Of course, the order of Vertex 1 of adding GDD is the same as order of  Vertex 2 of GDD $1$ of Row $6$. \\

 $\ \ \ \ \  \ \   
$\\

Remark: (a') is a GDD adding non classical GDD. GDDs from (a) to (v) are ones adding  classical GDDs.

{\rm (ii)}     Adding on Vertex 1 by a GDD in Table A1..\\

 $\ \ \ \ \  \ \   %
$\\  \\  by    Type 7, GDD $1$ of Row $16$ and GDD $2$ of Row $17$,  {\rm ord }   ($q_{11}) = 2$.\\

{\rm (iii)}   Cycle.

Every  quasi-affine cycle over GDD $1$ of Row $6$  consists of a GDD in Part {\rm (i)} and a GDD in Part {\rm (ii)}.  Of course, order of Vertex 3 of GDD in Part {\rm (i)} is the same as  order of Vertex 1 of GDD  Part {\rm (ii)}.
\\

$%
$\\ \\
by  GDD $3$ of Row $17$.

\subsection* {
GDD 2 of Row 6 in Table A1 }\label {sub3.2}

It is the same as GDD $1$ of Row $6$.

\subsection* {
GDD 1 of Row 8 in Table A1 }\label {sub3.2}

 {\rm (i)}      Adding on Vertex 2 by a GDD in Table A1.\\

 $\ \ \ \ \  \ \   %
$\\

{\rm (ii)}     Adding on Vertex 1 by a GDD in Table A1.\\

 $\ \ \ \ \  \ \   %
$\\ \\  It is empty.
 \\

$%
$\\ \\ $\lambda \in R_{12}$, $q ^{4}=1. $ \\ It is empty.\\ \\

$%
$\\ \\ $\lambda \in  R_{12},$ $q ^{4}=1. $\\ It is empty.\\

$%
$\\ \\  It is empty.\\

$%
$\\ \\   It is empty.\\

 $%
$\\ \\  $\lambda \in R_{12}$, $q ^{4}=1. $\\ It is empty.

$%
$\\ \\  $\lambda \in  R_{12},$ $q ^{4}=1. $ \\ It is empty.\\ \\

$%
$\\ \\  It is empty.\\

$%
$\\ \\  $q ^{4}=1. $ It is empty.\\

$%
$

\subsection* {
GDD 2 of Row 8 in Table A1 }\label {sub3.2}

{\rm (i)}    Adding on Vertex 2 by a GDD in Table A1.\\

 $\ \ \ \ \  \ \   %
$\\

{\rm (ii)}      Adding on Vertex 1 by a GDD in Table A1.\\

 $\ \ \ \ \  \ \   %
$

\subsection* {
GDD 3 of Row 8 in Table A1 }\label {sub3.2}

{\rm (i)}    Adding on Vertex 2 by a GDD in Table A1.
\\

 $\ \ \ \ \  \ \   %
$\\

{\rm (ii)}      Adding on Vertex 1 by a GDD in Table A1.\\

 $\ \ \ \ \  \ \   %
$

\subsection* {
GDD 4 of Row 8 in Table A1 }\label {sub3.2}

{\rm (i)}   Adding on Vertex 2 by a GDD in Table A1.
\\

 $\ \ \ \ \  \ \   %
$

\subsection* {
GDD 5 of Row 8 in Table A1 }\label {sub3.2}

{\rm (i)}    Adding on Vertex 2 by a GDD in Table A1.
\\

 $\ \ \ \ \  \ \   %
$

\subsection* {
GDD 1 of Row 9 in Table A1 }\label {sub3.2}

 {\rm (i)}     Adding on Vertex 2 by a GDD in Table A1.\\

 $\ \ \ \ \  \ \   %
$\\ \\ $q ^{4}=1. $  It is empty.
\\ \\

$%
$\\ \\  $q ^{4}=1. $ It is empty.\\

$%
$\\ \\ $q ^{4}=1. $ It is empty.\\

$%
$\\ \\ $q ^{4}=1. $ It is empty. \\

$%
$

\subsection* {
GDD 2 of Row 9 in Table A1 }\label {sub3.2}

{\rm (i)}    Adding on Vertex 2 by a GDD in Table A1.
\\

 $\ \ \ \ \  \ \   %
$

\subsection* {
GDD 3 of Row 9 in Table A1 }\label {sub3.2}

{\rm (i)}    Adding on Vertex 2 by a GDD in Table A1.
\\

 $\ \ \ \ \  \ \   %
$

\subsection* {
GDD 1 of Row 10 in Table A1 }\label {sub3.2}

 {\rm (i)}      Adding on Vertex 2 by a GDD in Table A1.\\

 $\ \ \ \ \  \ \   %
$

\subsection* {
GDD 2 of Row 10 in Table A1 }\label {sub3.2}

{\rm (i)}    Adding on Vertex 2 by a GDD in Table A1.
\\

 $\ \ \ \ \  \ \   %
$

\subsection* {
GDD 3 of Row 10 in Table A1 }\label {sub3.2}

{\rm (i)}    Adding on Vertex 2 by a GDD in Table A1.
\\

 $\ \ \ \ \  \ \   %
$

\subsection* {
GDD 1 of Row 11 in Table A1 }\label {sub3.2}

{\rm (i)}     Adding on Vertex 2 by a GDD in Table A1.\\

 $\ \ \ \ \  \ \   %
$\\ \\ $\mu \in R_{12}$, $q \in R_{18}$. \\  It is empty.\\

$%
$\\ \\ It is empty.
 \\

$%
$\\ \\ $q \in R_{18}$.
  It is empty.\\

$%
$\\ \\
  $q \in R_{6}$, $\lambda \in R_{12}$, $\mu \in R_{12}$, $\lambda \in R_{9}$,   $q \in  R_{18}$,\\  It is empty.\\

$%
$\\ \\
 $-\mu^{2} = q ^{3}, $ $-\mu ^{2} = \lambda^{3}, $   $-\lambda ^{}= q ^{}, $    $q \in R_{9}$, $\lambda \in R_{6}$, $\xi \in R_{3}$, $\mu \in R_{12}$,$q \in R_{18}$. \\ It is empty.\\

$%
$\\ \\
$\lambda \in R_{12}$, $\xi \in R_{3}$, $\mu \in R_{12}$, $\lambda \in  R_{3}$. \\
 It is empty.\\

$%
$\\ \\
 $\mu ^{3}= \mu ^{9}, $    $q ^{}= -1. $ \\
   It is empty.\\\\

$%
$

\subsection* {
GDD 1 of Row 12 in Table A1 }\label {sub3.2}

{\rm (i)}       Adding on Vertex 2 by a GDD in Table A1.\\

 $\ \ \ \ \  \ \   %
$

\subsection* {
GDD 2 of Row 12 in Table A1 }\label {sub3.2}

{\rm (i)}  Adding on Vertex 2 by a GDD in Table A1.
\\

 $\ \ \ \ \  \ \   %
$

\subsection* {
GDD 3 of Row 12 in Table A1 }\label {sub3.2}

{\rm (i)}   Adding on Vertex 2 by a GDD in Table A1.
\\

 $\ \ \ \ \  \ \   %
$

\subsection* {
GDD 1 of Row 13 in Table A1 }\label {sub3.2}

 {\rm (i)}    Adding on Vertex 2 by a GDD in Table A1.\\

 $\ \ \ \ \  \ \   %
$

\subsection* {
GDD 2 of Row 13 in Table A1 }\label {sub3.2}

{\ }

 $\ \ \ \ \  \ \   %
$

\subsection* {
GDD 3 of Row 13 in Table A1 }\label {sub3.2}

{\rm (i)}  Adding on Vertex 2 by a GDD in Table A1.
\\

 $\ \ \ \ \  \ \   %
$

\subsection* {
GDD 4 of Row 13 in Table A1 }\label {sub3.2}

{\rm (i)}   Adding on Vertex 2 by a GDD in Table A1.
\\

 $\ \ \ \ \  \ \   %
$

\subsection* {
GDD 1 of Row 14 in Table A1 }\label {sub3.2}

{\rm (i)}   Adding on Vertex 2 by a GDD in Table A1.
\\

 $\ \ \ \ \  \ \   %
$

\subsection* {
GDD 2 of Row 14 in Table A1 }\label {sub3, 2}

{\rm (i)}   Adding on Vertex 2 by a GDD in Table A1.
\\

 $\ \ \ \ \  \ \   %
$

\subsection* {
GDD 1 of Row 15 in Table A1 }\label {sub3, 2}

{\rm (i)}   Adding on Vertex 2 by a GDD in Table A1.
\\

 $\ \ \ \ \  \ \   %
$

\subsection* {
GDD 2 of Row 15 in Table A1 }\label {sub3, 2}

{\rm (i)}  Adding on Vertex 2 by a GDD in Table A1.\\

 $\ \ \ \ \  \ \   %
$

\subsection* {
GDD 3 of Row 15 in Table A1 }\label {sub3, 2}

{\rm (i)}  Adding on Vertex 2 by a GDD in Table A1.
\\

 $\ \ \ \ \  \ \   %
$

\subsection* {
GDD 4 of Row 15 in Table A1 }\label {sub3, 2}

{\rm (i)}  Adding on Vertex 2 by a GDD in Table A1.\\
\\

 $\ \ \ \ \  \ \   %
$

\subsection* {
GDD 1 of Row 16 in Table A1 }\label {sub3, 2}

  {\rm (i)}  Adding on Vertex 2 by a GDD in Table A1.\\

 $\ \ \ \ \  \ \   %
$\\ \\   It is empty. \\

  $%
$

\subsection* {
GDD 2 of Row 16 in Table A1 }\label {sub3, 2}

 {\rm (i)}  Adding on Vertex 2 by a GDD in Table A1.\\

 $\ \ \ \ \  \ \   %
$

\subsection* {
GDD 3 of Row 16 in Table A1 }\label {sub3, 2}
{\rm (i)}  Adding on Vertex 2 by a GDD in Table A1.
\\

 $\ \ \ \ \  \ \   %
$

\subsection* {
GDD 4 of Row 16 in Table A1 }\label {sub3, 2}

{\rm (i)}  Adding on Vertex 2 by a GDD in Table A1.
\\

 $\ \ \ \ \  \ \   %
$

 {\rm (ii)}   Adding on Vertex 1 by a GDD in Table A1.\\

 $\ \ \ \ \  \ \   %
$

\subsection* {
GDD 1 of Row 17 in Table A1 }\label {sub3, 2}

{\rm (i)}   Adding on Vertex 2 by a GDD in Table A1.
\\

 $\ \ \ \ \  \ \   %
$

\subsection* {
GDD 2 of Row 17 in Table A1 }\label {sub3, 2}

{\rm (i)}  Adding on Vertex 2 by a GDD in Table A1.\\

 $\ \ \ \ \  \ \   %
$

\subsection* {
Type 2 }\label {sub3, 2}

   $%
$\\

 {\rm (i)}  Adding on Vertex 2 by a GDD in Table A1. \\

 $\ \ \ \ \  \ \   %
$

\subsection* {
Type 2 with $-1$ }\label {sub3, 2}

   $%
$\\

 {\rm (i)} Adding on Vertex 2 by a GDD in Table A1.\\

 $\ \ \ \ \  \ \   %
$

\subsection* {
Type 4 }\label {sub3, 2}

   $%
$\\

 {\rm (i)} Adding on Vertex 2 by a GDD in Table A1.\\

 $\ \ \ \ \  \ \   %
$\\ \\ It is empty.

\subsection* {
Type 4 with $-1$ }\label {sub3, 2}

   $%
$\\

  {\rm (i)}   Adding on Vertex 2 by a GDD in Table A1.\\

 $\ \ \ \ \  \ \   %
$

\subsection* {
Type 7 }\label {sub3, 2}

 $%
$\\  is Type 7,   {\rm ord }   ($q_{33})=2$.
 \\

 {\rm (ii)}     Adding on Vertex 1 by a GDD in Table A1.\\

$%
$

\subsection* {
Type 7 with $-1$ }\label {sub3, 2}

   $%
$\\ is    Type 5 or  GDD $1$ of Row $9$.       {\rm ord }   ($q_{33})>1$. \\
   \\

 {\rm (ii)}     Adding on Vertex 1 by a GDD in Table A1.\\

 $\ \ \ \ \  \ \   %
$\\
It is  Type 7. \\

$%
$\\ \\
It is  GDD $8$ of Row $17$. \\

$%
$

\subsection* {
Type 7 with two $-1$}\label {sub3, 2}

  $%
$\\ \\  Case 2: $ q^{-1} =\mu ^{2}$, $ \lambda^{-1} =\mu^{-1}$.

\section {    Quasi-affine  Nichols  (Lie braided) Algebras with Rank $2$    } \label {s3}

In this section all quasi-affine  connected Generalized Dynkin Diagram with rank $= 2$ are found,  All quasi-affine  Nichols (Lie braided) algebras with rank $ 2$  are also found.

({\rm ord }   ($q_{11})$,  {\rm ord }   ($\widetilde{q}_{12})$,  {\rm ord }   ($q_{22})$)
  is called the order of GDD
$%
$\\

Let $\Omega  := \{   x \mid x \hbox { is  order of   an  arithmetic   GDD which is not classical.   }  GDD $ 1$
  \hbox { of Row } $6$ \\  \hbox { and  GDD } $1$   \hbox { of Row } $11$ \}$.
  In other words, $\Omega  = \{   x \mid x \hbox { is order of   an  arithmetic   GDD } \hbox { in  Row } $ 8, $  \hbox { Row} $ 9, $      \hbox { Row} $ 10, $  \hbox {   Row } $ 12, $  \hbox { Row} $ 13, $  \hbox { Row} $ 14, $   \hbox { Row} $ 15, $ \hbox { Row} $ 16, $  \hbox { Row} $ 17 $   \}$.

\begin  {Proposition} \label {3.1}
        A connected  GDD with rank 2 is quasi-affine  if and only if it is not   an arithmetic GDD.   \end {Proposition}

 \begin {Lemma} \label {3.2}  If   there exist two places  where simple chain conditions do not hold, then the GDD is quasi-affine
except
  GDD $1$ of Row $8$,  GDD $1$ of Row $9$ , GDD $1$ of Row $10$, GDD $1$ of Row $13$,   GDD $1$ of Row $16. $
\end {Lemma}

 \begin {Theorem} \label {Mainrank2}
        A connected  GDD with rank $= 2$ is quasi-affine   if and only if it is  one of following lists.

  {\rm (i)}  Assume that  order of the GDD  belongs to $\Omega$.

  (1)  GDD with  order  $(3.4.3)$ except
   GDD $1$ of Row $8$: $\begin{picture}(100,      15)  \put(-20,      -1){ }

\put(60,      1){\makebox(0,      0)[t]{$\bullet$}}

\put(28,      -1){\line(1,      0){33}}
\put(27,      1){\makebox(0,     0)[t]{$\bullet$}}



\put(22,     10){$-q ^{-2}$}
\put(38,      5){$-q^{3}$}

\put(63,      10){$-q^{2}$}

  \ \ \ \ \ \ \ \ \ \ \ \ \ \ \ \ \ \ \ {$q\in R_{12}$.  }
\end{picture}$\\

 (2)  GDD with order  $( 3, 12, 2)$ except   GDD $2$ of Row $8$:
$\begin{picture}(100,      15)  \put(-20,      -1){}

\put(60,      1){\makebox(0,      0)[t]{$\bullet$}}

\put(28,      -1){\line(1,      0){33}}
\put(27,      1){\makebox(0,     0)[t]{$\bullet$}}



\put(18,     10){$-q ^{-2}$}
\put(38,      5){$q^{-1}$}

\put(63,      10){$-1$}

  \ \ \ \ \ \ \ \ \ \ \ \ \ \ \ \ \ \ \ {$q\in R_{12}$. }
\end{picture}$\\

 and  GDD $3$ of Row $8$:   $\begin{picture}(100,      15)  \put(-20,      -1){}

\put(60,      1){\makebox(0,      0)[t]{$\bullet$}}

\put(28,      -1){\line(1,      0){33}}
\put(27,      1){\makebox(0,     0)[t]{$\bullet$}}



\put(22,     10){$-q ^{2}$}
\put(38,      5){-$q^{}$}

\put(63,      10){$-1$}

  \ \ \ \ \ \ \ \ \ \ \ \ \ \ \ \ \ \ \ {$q\in R_{12}$.  }
\end{picture}$\\

 (3)  GDD with  order   $( 4, 12, 2)$ except  GDD $4$ of Row $8$:
$\begin{picture}(100,      15)  \put(-20,      -1){ }

\put(60,      1){\makebox(0,      0)[t]{$\bullet$}}

\put(28,      -1){\line(1,      0){33}}
\put(27,      1){\makebox(0,     0)[t]{$\bullet$}}



\put(18,     10){$-q ^{2}$}
\put(38,      5){$q^{}$}

\put(63,      10){$-1$}

  \ \ \ \ \ \ \ \ \ \ \ \ \ \ \ \ \ \ \ {$q\in R_{12}$. }
\end{picture}$\\

and GDD $5$ of Row $8$:
$\begin{picture}(100,      15)  \put(-20,      -1){ }

\put(60,      1){\makebox(0,      0)[t]{$\bullet$}}

\put(28,      -1){\line(1,      0){33}}
\put(27,      1){\makebox(0,     0)[t]{$\bullet$}}



\put(22,     10){$-q ^{3}$}
\put(38,      5){$-q^{-1}$}

\put(63,      10){$-1$}

  \ \ \ \ \ \ \ \ \ \ \ \ \ \ \ \ \ \ \ {$q\in R_{12}$.  }
\end{picture}$\\

 (4)  GDD with  order    GDD is   $( 3, 12, 3)$ except  GDD $1$ of Row $9$:   $\begin{picture}(100,      15)  \put(-20,      -1){ }

\put(60,      1){\makebox(0,      0)[t]{$\bullet$}}

\put(28,      -1){\line(1,      0){33}}
\put(27,      1){\makebox(0,     0)[t]{$\bullet$}}



\put(22,     10){$-q ^{2}$}
\put(38,      5){$q^{}$}

\put(60,      10){$-q^{2}$}

  \ \ \ \ \ \ \ \ \ \ \ \ \ \ \ \ \ \ \ {$q\in R_{12}$.  }
\end{picture}$\\

 (5)  GDD with  order     $( 3, 4, 2)$ except  GDD $2$ of Row $9$:   $\begin{picture}(100,      15)  \put(-20,      -1){ }

\put(60,      1){\makebox(0,      0)[t]{$\bullet$}}

\put(28,      -1){\line(1,      0){33}}
\put(27,      1){\makebox(0,     0)[t]{$\bullet$}}



\put(22,     10){$-q ^{2}$}
\put(38,      5){$q^{3}$}

\put(60,      10){$-1$}

  \ \ \ \ \ \ \ \ \ \ \ \ \ \ \ \ \ \ \ {$q\in R_{12}$.  }
\end{picture}$\\

 (6)  GDD with  order    $( 12, 4, 2)$ except  GDD $3$ of Row $9$:   $\begin{picture}(100,      15)  \put(-20,      -1){ }

\put(60,      1){\makebox(0,      0)[t]{$\bullet$}}

\put(28,      -1){\line(1,      0){33}}
\put(27,      1){\makebox(0,     0)[t]{$\bullet$}}



\put(22,     10){$-q ^{-1}$}
\put(38,      5){$-q^{3}$}

\put(60,      10){$-1$}

  \ \ \ \ \ \ \ \ \ \ \ \ \ \ \ \ \ \ \ {$q\in R_{12}$.  }
\end{picture}$\\

 (7)  GDD with  order    $( 18, 9, 3)$ except  GDD $1$ of Row $10$:   $\begin{picture}(100,      15)  \put(-20,      -1){ }

\put(60,      1){\makebox(0,      0)[t]{$\bullet$}}

\put(28,      -1){\line(1,      0){33}}
\put(27,      1){\makebox(0,     0)[t]{$\bullet$}}



\put(22,     10){$-q ^{}$}
\put(38,      5){$q^{-2}$}

\put(60,      10){$q^{3}$}

  \ \ \ \ \ \ \ \ \ \ \ \ \ \ \ \ \ \ \ {$q\in R_{9}$.  }
\end{picture}$\\

 (8)  GDD with  order    $( 3, 9, 2)$ except  GDD $2$ of Row $10$:   $\begin{picture}(100,      15)  \put(-20,      -1){}

\put(60,      1){\makebox(0,      0)[t]{$\bullet$}}

\put(28,      -1){\line(1,      0){33}}
\put(27,      1){\makebox(0,     0)[t]{$\bullet$}}



\put(22,     10){$q ^{3}$}
\put(38,      5){$q^{-1}$}

\put(60,      10){$-1$}

  \ \ \ \ \ \ \ \ \ \ \ \ \ \ \ \ \ \ \ {$q\in R_{9}$.  }
\end{picture}$\\

 (9)  GDD with  order   $( 18, 9, 2)$ except GDD $3$ of Row $10$:   $\begin{picture}(100,      15)  \put(-20,      -1){ }

\put(60,      1){\makebox(0,      0)[t]{$\bullet$}}

\put(28,      -1){\line(1,      0){33}}
\put(27,      1){\makebox(0,     0)[t]{$\bullet$}}



\put(22,     10){$-q ^{2}$}
\put(38,      5){$q^{}$}

\put(60,      10){$-1$}

  \ \ \ \ \ \ \ \ \ \ \ \ \ \ \ \ \ \ \ {$q\in R_{9}$.  }
\end{picture}$\\

 (10)  GDD with  order     $( 4, 8, 8)$  except  GDD $1$ of Row $12$:   $\begin{picture}(100,      15)  \put(-20,      -1){ }

\put(60,      1){\makebox(0,      0)[t]{$\bullet$}}

\put(28,      -1){\line(1,      0){33}}
\put(27,      1){\makebox(0,     0)[t]{$\bullet$}}



\put(22,     10){$q ^{2}$}
\put(38,      5){$q^{}$}

\put(60,      10){$q ^{-1}$}

  \ \ \ \ \ \ \ \ \ \ \ \ \ \ \ \ \ \ \ {$q\in R_{8}$.  }
\end{picture}$\\

 (11)  GDD with  order   $( 4, 8, 2)$ except  GDD $2$ of Row $12$:   $\begin{picture}(100,      15)  \put(-20,      -1){}

\put(60,      1){\makebox(0,      0)[t]{$\bullet$}}

\put(28,      -1){\line(1,      0){33}}
\put(27,      1){\makebox(0,     0)[t]{$\bullet$}}



\put(22,     10){$q ^{2}$}
\put(38,      5){$-q^{-1}$}

\put(60,      10){${-1}$}

  \ \ \ \ \ \ \ \ \ \ \ \ \ \ \ \ \ \ \ {$q\in R_{8}$.  }
\end{picture}$\\

 (12)  GDD with  order    $( 8, 8, 2)$ except  GDD $3$ of Row $12$:   $\begin{picture}(100,      15)  \put(-20,      -1){ }

\put(60,      1){\makebox(0,      0)[t]{$\bullet$}}

\put(28,      -1){\line(1,      0){33}}
\put(27,      1){\makebox(0,     0)[t]{$\bullet$}}



\put(22,     10){$q ^{}$}
\put(38,      5){$-q^{}$}

\put(60,      10){${-1}$}

  \ \ \ \ \ \ \ \ \ \ \ \ \ \ \ \ \ \ \ {$q\in R_{8}$.  }
\end{picture}$\\

  (13)  GDD with  order     $( 4, 24, 3)$  except  GDD $1$ of Row $13$:   $\begin{picture}(100,      15)  \put(-20,      -1){ }

\put(60,      1){\makebox(0,      0)[t]{$\bullet$}}

\put(28,      -1){\line(1,      0){33}}
\put(27,      1){\makebox(0,     0)[t]{$\bullet$}}



\put(22,     10){$q ^{6}$}
\put(38,      5){$-q^{-1}$}

\put(60,      10){$-q ^{-4}$}

  \ \ \ \ \ \ \ \ \ \ \ \ \ \ \ \ \ \ \ {$q\in R_{24}$.  }
\end{picture}$\\

 (14)  GDD with  order   $( 4, 24, 24)$ except  GDD $2$ of Row $13$:   $\begin{picture}(100,      15)  \put(-20,      -1){ }

\put(60,      1){\makebox(0,      0)[t]{$\bullet$}}

\put(28,      -1){\line(1,      0){33}}
\put(27,      1){\makebox(0,     0)[t]{$\bullet$}}



\put(22,     10){$q ^{6}$}
\put(38,      5){$q^{}$}

\put(60,      10){$q ^{-1}$}

  \ \ \ \ \ \ \ \ \ \ \ \ \ \ \ \ \ \ \ {$q\in R_{24}$.  }
\end{picture}$\\

 (15)  GDD with  order     $( 3, 24, 2)$ except  GDD $3$ of Row $13$:   $\begin{picture}(100,      15)  \put(-20,      -1){ }

\put(60,      1){\makebox(0,      0)[t]{$\bullet$}}

\put(28,      -1){\line(1,      0){33}}
\put(27,      1){\makebox(0,     0)[t]{$\bullet$}}



\put(18,     10){$-q ^{-4}$}
\put(42,      5){$q^{5}$}

\put(60,      10){$-1$}

  \ \ \ \ \ \ \ \ \ \ \ \ \ \ \ \ \ \ \ {$q\in R_{24}$.  }
\end{picture}$\\

 (16)  GDD with  order   $( 24, 24, 2)$ except  GDD $4$ of Row $13$:   $\begin{picture}(100,      15)  \put(-20,      -1){ }

\put(60,      1){\makebox(0,      0)[t]{$\bullet$}}

\put(28,      -1){\line(1,      0){33}}
\put(27,      1){\makebox(0,     0)[t]{$\bullet$}}



\put(22,     10){$q ^{}$}
\put(44,      5){$q^{-5}$}

\put(60,      10){$-1$}

  \ \ \ \ \ \ \ \ \ \ \ \ \ \ \ \ \ \ \ {$q\in R_{24}$.  }
\end{picture}$\\

 (17)  GDD with  order     $( 5, 5, 2)$ except  GDD $1$ of Row $14$:   $\begin{picture}(100,      15)  \put(-20,      -1){}

\put(60,      1){\makebox(0,      0)[t]{$\bullet$}}

\put(28,      -1){\line(1,      0){33}}
\put(27,      1){\makebox(0,     0)[t]{$\bullet$}}



\put(22,     10){$q ^{}$}
\put(38,      5){$q^{2}$}

\put(60,      10){$-1$}

  \ \ \ \ \ \ \ \ \ \ \ \ \ \ \ \ \ \ \ {$q\in R_{5}$.  }
\end{picture}$ \\

 (18)  GDD with  order   $( 10, 5, 2)$ except  GDD $2$ of Row $14$:
   $\begin{picture}(100,      15)  \put(-20,      -1){ }

\put(60,      1){\makebox(0,      0)[t]{$\bullet$}}

\put(28,      -1){\line(1,      0){33}}
\put(27,      1){\makebox(0,     0)[t]{$\bullet$}}



\put(14,     10){$-q ^{-2}$}
\put(38,      5){$q^{-2}$}

\put(60,      10){$-1$}

  \ \ \ \ \ \ \ \ \ \ \ \ \ \ \ \ \ \ \ {$q\in R_{5}$.  }
\end{picture}$\\

 (19)  GDD with  order    $( 20, 20, 2)$ except  GDD $1$ of Row $15$:   $\begin{picture}(100,      15)  \put(-20,      -1){ }

\put(60,      1){\makebox(0,      0)[t]{$\bullet$}}

\put(28,      -1){\line(1,      0){33}}
\put(27,      1){\makebox(0,     0)[t]{$\bullet$}}



\put(22,     10){$q ^{}$}
\put(44,      5){$q^{-3}$}

\put(60,      10){$-1$}

  \ \ \ \ \ \ \ \ \ \ \ \ \ \ \ \ \ \ \ {$q\in R_{20}$.  }
\end{picture}$\\

and  GDD $2$ of Row $15$:   $\begin{picture}(100,      15)  \put(-20,      -1){ }

\put(60,      1){\makebox(0,      0)[t]{$\bullet$}}

\put(28,      -1){\line(1,      0){33}}
\put(27,      1){\makebox(0,     0)[t]{$\bullet$}}



\put(22,     10){$-q ^{}$}
\put(38,      5){$-q^{-3}$}

\put(60,      10){$-1$}

  \ \ \ \ \ \ \ \ \ \ \ \ \ \ \ \ \ \ \ {$q\in R_{20}$.  }
\end{picture}$\\

 (20)  GDD with  order    $( 5, 20, 2)$ except  GDD $3$ of Row $15$:   $\begin{picture}(100,      15)  \put(-20,      -1){ }

\put(60,      1){\makebox(0,      0)[t]{$\bullet$}}

\put(28,      -1){\line(1,      0){33}}
\put(27,      1){\makebox(0,     0)[t]{$\bullet$}}



\put(18,     10){$-q ^{-2}$}
\put(38,      5){$q^{3}$}

\put(60,      10){$-1$}

  \ \ \ \ \ \ \ \ \ \ \ \ \ \ \ \ \ \ \ {$q\in R_{20}$. }
\end{picture}$\\

 and  GDD $4$ of Row $15$:   $\begin{picture}(100,      15)  \put(-20,      -1){ }

\put(60,      1){\makebox(0,      0)[t]{$\bullet$}}

\put(28,      -1){\line(1,      0){33}}
\put(27,      1){\makebox(0,     0)[t]{$\bullet$}}



\put(22,     10){$-q ^{-2}$}
\put(38,      5){$-q^{3}$}

\put(60,      10){$-1$}

  \ \ \ \ \ \ \ \ \ \ \ \ \ \ \ \ \ \ \ {$q\in R_{20}$.  }
\end{picture}$\\

  (21)  GDD with  order     $( 30, 10, 3)$ except  GDD $1$ of Row $16$:   $\begin{picture}(100,      15)  \put(-20,      -1){}

\put(60,      1){\makebox(0,      0)[t]{$\bullet$}}

\put(28,      -1){\line(1,      0){33}}
\put(27,      1){\makebox(0,     0)[t]{$\bullet$}}



\put(12,     10){$-q ^{}$}
\put(28,      5){$-q^{-3}$}

\put(60,      10){$q^{5}$}

  \ \ \ \ \ \ \ \ \ \ \ \ \ \ \ \ \ \ \ {$q\in R_{15}$.  }
\end{picture}$\\

 (22)  GDD with  order    $( 5, 30, 30)$ except   GDD $2$ of Row $16$:   $\begin{picture}(100,      15)  \put(-20,      -1){ }

\put(60,      1){\makebox(0,      0)[t]{$\bullet$}}

\put(28,      -1){\line(1,      0){33}}
\put(27,      1){\makebox(0,     0)[t]{$\bullet$}}



\put(22,     10){$q ^{3}$}
\put(38,      5){$-q^{4}$}

\put(60,      10){$-q^{-4}$}

  \ \ \ \ \ \ \ \ \ \ \ \ \ \ \ \ \ \ \ {$q\in R_{15}$.  }
\end{picture}$\\

  (23)  GDD with  order    $( 3, 30, 2)$ except  GDD $3$ of Row $16$:   $\begin{picture}(100,      15)  \put(-20,      -1){}

\put(60,      1){\makebox(0,      0)[t]{$\bullet$}}

\put(28,      -1){\line(1,      0){33}}
\put(27,      1){\makebox(0,     0)[t]{$\bullet$}}



\put(22,     10){$q ^{5}$}
\put(34,      5){$-q^{-2}$}

\put(60,      10){$-1$}

  \ \ \ \ \ \ \ \ \ \ \ \ \ \ \ \ \ \ \ {$q\in R_{15}$. }
\end{picture}$\\

(24)  GDD with  order    $( 5, 30, 2)$ except  GDD $4$ of Row $16$:   $\begin{picture}(100,      15)  \put(-20,      -1){ }

\put(60,      1){\makebox(0,      0)[t]{$\bullet$}}

\put(28,      -1){\line(1,      0){33}}
\put(27,      1){\makebox(0,     0)[t]{$\bullet$}}



\put(22,     10){$q ^{3}$}
\put(38,      5){$-q^{2}$}

\put(60,      10){$-1$}

  \ \ \ \ \ \ \ \ \ \ \ \ \ \ \ \ \ \ \ {$q\in R_{15}$.  }
\end{picture}$\\

    (24)  GDD with  order    $( 14, 14, 2)$ except  GDD $1$ of Row $17$:   $\begin{picture}(100,      15)  \put(-20,      -1){ }

\put(60,      1){\makebox(0,      0)[t]{$\bullet$}}

\put(28,      -1){\line(1,      0){33}}
\put(27,      1){\makebox(0,     0)[t]{$\bullet$}}



\put(22,     10){$-q ^{}$}
\put(38,      5){$-q^{-3}$}

\put(60,      10){$-1$}

  \ \ \ \ \ \ \ \ \ \ \ \ \ \ \ \ \ \ \ {$q\in R_{7}$.  }
\end{picture}$\\

 and     GDD $2$ of Row $17$:   $\begin{picture}(100,      15)  \put(-20,      -1){}

\put(60,      1){\makebox(0,      0)[t]{$\bullet$}}

\put(28,      -1){\line(1,      0){33}}
\put(27,      1){\makebox(0,     0)[t]{$\bullet$}}



\put(22,     10){$-q ^{-2}$}
\put(38,      5){$-q^{3}$}

\put(60,      10){$-1$}

  \ \ \ \ \ \ \ \ \ \ \ \ \ \ \ \ \ \ \ {$q\in R_{7}$.  }
\end{picture}$\\

 {\rm (ii)}   GDD which    order  does not belong to $\Omega$ and in which there exist two places  where simple chain conditions do not hold.

{\rm (iii)}    GDD which    order  does not belong to $\Omega$ and in which  there exists only one place  where simple chain condition does not hold with $q_{22}  \widetilde{q}_{12} =1  $  except Case 1: $q_{11}^{-2} =  \widetilde{q}_{12},   $ $q_{11}^2 \not= 1$;  Case 2: $q_{11}^{} =  -\widetilde{q}_{12},   $
 $q_{11} \in  R_{3}$; Case 3: $q_{11}^{-3} =  \widetilde{q}_{12},    $
 {\rm ord }   ($q_{11} ) > 3$;  Case 4: $q_{11} \in R_{3}$,
 $q_{22} \in  R_2$  or {\rm ord }   ($q_{22} ) > 3$.

\end {Theorem}

{\bf Proof.  } It follows from  Lemma \ref {3.2}.   \hfill $\Box$

{Remark:}  If the order of  GDD does not belong to $\Omega$ and  there exists only one place  where simple chain condition does not hold, then the GDD is  quasi-affine  if and only if GDD is not  any one of following lists: Type 2, Type 4, GDD $1$ of Row $6$, GDD $1$ of Row $11$.

In other words,  {\rm (iii)}   in Theorem \ref {Mainrank2} can be denoted:   GDD which    order  does not belong to $\Omega$ and in which  there exists only one place  where simple chain condition does not hold except   Type 2, Type 4, GDD $1$ of Row $6$ and GDD $1$ of Row $11$.

\section {Appendix}

\subsection* {Omitting one vertex in  GDDs with rank 3}

Omitting 1 in  GDD $1$ of Row $7$ is  GDD $1$ of Row $11$.   Omitting 3 in  GDD $1$ of Row $7$ is  Type 7.

Omitting 1 in  GDD $2$ of Row $7$ is  Type 7.   Omitting 3 in  GDD $2$ of Row $7$ is  Type 7.

Omitting 1 in  GDD $4$ of Row $7$ is  Type 3.   Omitting 3 in  GDD $4$ of Row $7$ is  Type 7.

Omitting 1 in  GDD $1$ of Row $9$ is  Type 7.   Omitting 3 in  GDD $1$ of Row $9$ is  Type 7.

Omitting 1 in  GDD $1$ of Row $13$ is  Type 2.   Omitting 3 in  GDD $1$ of Row $13$ is  Type 7.

Omitting 1 in  GDD $2$ of Row $13$ is  Type 3.   Omitting 3 in  GDD $2$ of Row $13$ is  Type 3 when $q \in  R_{3}$.     Omitting 3 in  GDD $2$ of Row $13$ is  Type 4 when $q \in  R_{6}$.

Omitting 1 in  GDD $1$ of Row $15$ is  Type 3.   Omitting 3 in  GDD $1$ of Row $15$ is  Type 7.

Omitting 1 in  GDD $2$ of Row $15$ is  Type 7.   Omitting 3 in  GDD $2$ of Row $15$ is  Type 7.

Omitting 1 in  GDD $4$ of Row $15$ is  Type 3.   Omitting 3 in  GDD $4$ of Row $15$ is  Type 3.

Omitting 1 in  GDD $1$ of Row $16$ is   GDD $1$ of Row $6$ with  $q= \zeta $ .   Omitting 3 in  GDD $1$ of Row $16$ is  Type 7.

Omitting 1 in  GDD $2$ of Row $16$ is   Type 7.   Omitting 3 in  GDD $2$ of Row $16$ is  Type 7.

Omitting 1 in  GDD $4$ of Row $16$ is   Type 7.   Omitting 3 in  GDD $4$ of Row $16$ is  GDD $1$ of Row $6$ with  $q= -1 $.   $\xi= \zeta $.

Omitting 1 in  GDD $5$ of Row $16$ is   Type 7.   Omitting 3 in  GDD $5$ of Row $16$ is  GDD $1$ of Row $6$ with  $q= -\zeta $.   $\xi= \zeta $.

Omitting 1 in  GDD $1$ of Row $17$ is   Type 7.   Omitting 3 in  GDD $1$ of Row $17$ is  Type 7.

Omitting 1 in  GDD $2$ of Row $17$ is   Type 7.   Omitting 3 in  GDD $2$ of Row $17$ is  GDD $1$ of Row $6$ with  $q= -1 $.   $\xi= \zeta $.

Omitting 1 in  GDD $4$ of Row $17$ is   Type 7.   Omitting 3 in  GDD $4$ of Row $17$ is  Type 7.

Omitting 1 in  GDD $5$ of Row $17$ is   Type 7.   Omitting 3 in  GDD $5$ of Row $17$ is  Type 3.

Omitting 1 in  GDD $6$ of Row $17$ is   Type 4.   Omitting 3 in  GDD $6$ of Row $17$ is  Type 2.

Omitting 1 in  GDD $7$ of Row $17$ is   Type 4.   Omitting 3 in  GDD $6$ of Row $17$ is  Type 7.

Omitting 1 in  GDD $9$ of Row $17$ is   Type 4.   Omitting 3 in  GDD $9$ of Row $17$ is  Type 7.

Omitting 1 in  GDD $1$ of Row $18$ is  GDD $1$ of Row $6$ with  $q= \zeta $.   $\xi= \zeta^{-3} $ .   Omitting 3 in  GDD $1$ of Row $18$ is  Type 7.

Omitting 1 in  GDD $2$ of Row $18$ is  GDD $1$ of Row $6$ with  $q= \zeta ^{-4} $.   $\xi= \zeta^{-3} $ .   Omitting 3 in  GDD $2$ of Row $18$ is  Type 2.

For Cycle.

Omitting 2 in  GDD $3$ of Row $7$ is   Type 2.

Omitting 2 in  GDD $3$ of Row $16$ is  GDD $1$ of Row $6$ with  $q= -1.  $ $ \xi= \zeta ^{} $.

Omitting 2 in  GDD $8$ of Row $17$ is   Type 4.

Omitting 3 in  GDD $3$ of Row $17$ is    GDD $1$ of Row $6$ with  $q= -\zeta.  $ $ \xi= \zeta ^{} $.


\begin{thebibliography}{BD99}

\bibitem [AS02]{AS02}
 N. Andruskiewitsch,  H.-J. Schneider,   Pointed Hopf algebras,  in New Directions in Hopf Algebras,  MSRI Series, pp.1-68,  Cambridge Univ. Press, 2002.








\bibitem[AS10] {AS10} N. Andruskiewitsch,   H.-J. Schneider,    On the classification of finite-dimensional pointed Hopf algebras,
 Ann. Math.   {\bf 171}  (2010),   375-417.
\bibitem[AS00] {AS00} N. Andruskiewitsch and H.J. Schneider,
Finite quantum groups and Cartan matrices,    Adv. Math. {\bf 154}  (2000),    1-45.

\bibitem [AHS08] {AHS08} N. Andruskiewitsch,      I. Heckenberger and   H.J. Schneider,
  The Nichols algebra of a semisimple Yetter-Drinfeld module,
 Amer. J. Math. {\bf 132}  (2010),     1493-1547.



\bibitem[Ha] {Ha} F. Harary, Graph Theory,
Addison-Wesley,  USA,  1969.

\bibitem[He05] {He05} I. Heckenberger,      Nichols algebras of diagonal type and arithmetic root systems,   Habilitation thesis, Leipzig, 2005.

\bibitem [He06a] {He06a}I. Heckenberger,  Classification of arithmetic
root systems,  Adv. Math.  {\bf 220} (2009),  59-124.

\bibitem [He06b]{He06b} I. Heckenberger,  The Weyl-Brandt groupoid of a Nichols algebra
of diagonal  Type,  Invent. Math. {\bf 164} (2006),  175--188.






\bibitem  [Ka90] {Ka90}V. G. Kac. Infinite dimensional Lie Algebras. London:Cambridge Univ, 1990.


\bibitem [Sc79]{Sc79} M. Scheunert, Generalized Lie algebras, J. Math. Phys. {\bf 20} (1979), 712-720.



\bibitem [WZZ15a]{WZZ15a} W. Wu,   S. Zhang and   Y.-Z. Zhang,    Relationship between Nichols braided Lie
algebras and Nichols algebras,  J. Lie Theory {\bf 25} (2015),     45-63.

\bibitem [WZZ15b]{WZZ15b} W. Wu,   S. Zhang and   Y.-Z. Zhang,    On  Nichols (braided) Lie algebras,  Int. J. Math. {\bf 26} (2015),  1550082. 


\end{thebibliography}
\end {document}